\documentclass[preprint]{elsarticle}
\pdfoutput=1

\usepackage{lineno,hyperref}
\modulolinenumbers[1]
\usepackage{float}
\usepackage{graphicx}
\usepackage{mathrsfs}
\usepackage{graphics}
\usepackage{float}
\usepackage{amsmath}
\usepackage{amsfonts}
\usepackage{optidef}
\usepackage[subnum]{cases}
\usepackage{svg}
\usepackage{enumitem}

\usepackage{threeparttable}

\usepackage{optidef}
\usepackage{graphicx}
\usepackage{float}
\usepackage{graphicx}
\usepackage{mathrsfs}
\usepackage{float}
\usepackage{amsmath}
\usepackage{amsfonts}
\usepackage{booktabs}
\usepackage{algorithm}
\usepackage[noend]{algpseudocode}
\usepackage{multirow}
\algdef{SE}[DOWHILE]{Do}{doWhile}{\algorithmicdo}[1]{\algorithmicwhile\ #1}

\usepackage[letterpaper, portrait, margin=1in]{geometry}

\usepackage{optidef}

\usepackage{amsthm}

\theoremstyle{definition}

\makeatletter
\newcommand*{\centerfloat}{%
  \parindent \z@
  \leftskip \z@ \@plus 1fil \@minus \marginparwidth
  \rightskip \leftskip
  \parfillskip \z@skip}
\makeatother

\usepackage{setspace}
\usepackage{natbib}
\usepackage{algorithm}
\usepackage[noend]{algpseudocode}
\algdef{SE}[DOWHILE]{Do}{doWhile}{\algorithmicdo}[1]{\algorithmicwhile\ #1}

\usepackage{appendix}

\usepackage{array}
\usepackage{threeparttable}
\usepackage{multirow}
\usepackage{rotating}
\usepackage{makecell}

\usepackage{comment}

\newsavebox{\measurebox}

\journal{}

\biboptions{authoryear}

\usepackage{makecell}
\usepackage{longtable}
\usepackage{caption}
\usepackage{subcaption}
\usepackage{multirow}

\begin{document}

\begin{frontmatter}


\title{%
       Understanding and Analyzing the Influential Factors \\ 
       on Relocation of Shared Bikes}

\author[1]{Xinling Li}
\author[1]{Yu Shen\corref{mycorrespondingauthor}}\ead{yshen@tongji.edu.cn}
\author[1]{Chi Xie\corref{mycorrespondingauthor}}\ead{chi.xie@tongji.edu.cn}
\author[2]{Xiaohu Zhang}
\author[1]{Hanjun Fu}

\address[1]{Key Laboratory of Road and Traffic Engineering of the Ministry of Education, Tongji University, Shanghai 201804, China}
\address[2]{Senseable City Laboratory, Massachusetts Institute of Technology, Cambridge, MA 02139}
\cortext[mycorrespondingauthor]{Corresponding authors}

\begin{abstract}
Dockless sharing bike programs represent an important step towards smarter mobility service and boost transportation service accessibility. To enhance the service quality, bike fleet relocation is widely applied to achieve a better bike-rider balance across different areas. The primary challenge lies in devising a systematic approach that accommodates the dynamic nature of the demand during the day. In this study, we introduce a network flow model designed to address the optimal relocation problem of shared bikes. The model is tested with the dockless shared bike usage data from Yishun, Singapore, and demonstrated to be effective in providing resonable bike relocation scheme. Based on the test result, a series of sensitivity analyses were performed to investigate the impact of the relocation cost, the number of bikes and truck trikes, and the usage price on bike relocation. The sensitivity analysis reveals the connection between the profitability of the system and the analyzed factors. This work not only offers a modeling framework to initiate and manage a bikesharing service but also provides valuable insights for determining the number of bikes and trikes as well as price schemes. Additionally, some regulatory policies for the effective operation of bikesharing systems are suggested.
\end{abstract}

\begin{keyword}
(Dockless) bikesharing \sep Fleet Relocation \sep Network Flow Model


\end{keyword}

\end{frontmatter}


\section{Introduction and Background}\label{intro}

Over the past few decades, the proliferation of bikesharing programs (BSPs) has played a vital role in reducing vehicle emissions and enhancing the first-/last-mile connectivity. The evolution of bikesharing systems can be categorized into three distinct generations: 1) the unlocked free-of-charge free-floating shared bikes, 2) the coin deposit-based bikesharing with docking stations, and 3) the information and communication technology (ICT)-based bikesharing system with docking stations \citep{demaio2009bike,shaheen2010bikesharing,fishman2016bikeshare}. The first and second generations, while pioneering, were fundamentally basic. They lacked an effective information-gathering device to help manage fleets and prevent theft and, therefore, ultimately leading to their demise. The third generation marked a significant advancement, introducing docking stations and automated card payments for better management of bike fleets. ICT was also integrated to enable precise bike tracking. Along with advances in technology, these programs also get upgraded to provide better services. In recent years, with the development of mobile technology, a new generation of bikesharing systems, the dockless bikesharing system, which is equipped with mobile payment modules and GPS chips, emerged. This generation widely spread worldwide and became a new trend and a main part of BSPs.

Dockless bikesharing systems provide detailed information on the condition of all shared bikes, seemingly simplifying their operation and maintenance. However, in practice, dockless BSPs still hold a considerable risk of failure, such as in the cases of \textit{ofo} in China and \textit{oBike} in Singapore. One contributing factor to such failure is the substantial influx of shared bikes from diverse operators, demanding a nuanced understanding of the demand patterns. An oversupply of shared bikes exacerbates the issue by consuming public space, while an undersupply diminishes the operators’ market share and profitability, ultimately undermining the sustainability of dockless bikesharing services. In the following section, we will discuss the case of bike oversupply by utilizing data from one of the largest bikesharing operators in Singapore to support the development of our model.

\subsection{Dockless Bikesharing in Singapore}\label{dbs_sg}
\begin{figure}[ht]
\centering
\includegraphics[width= 0.9 \textwidth]{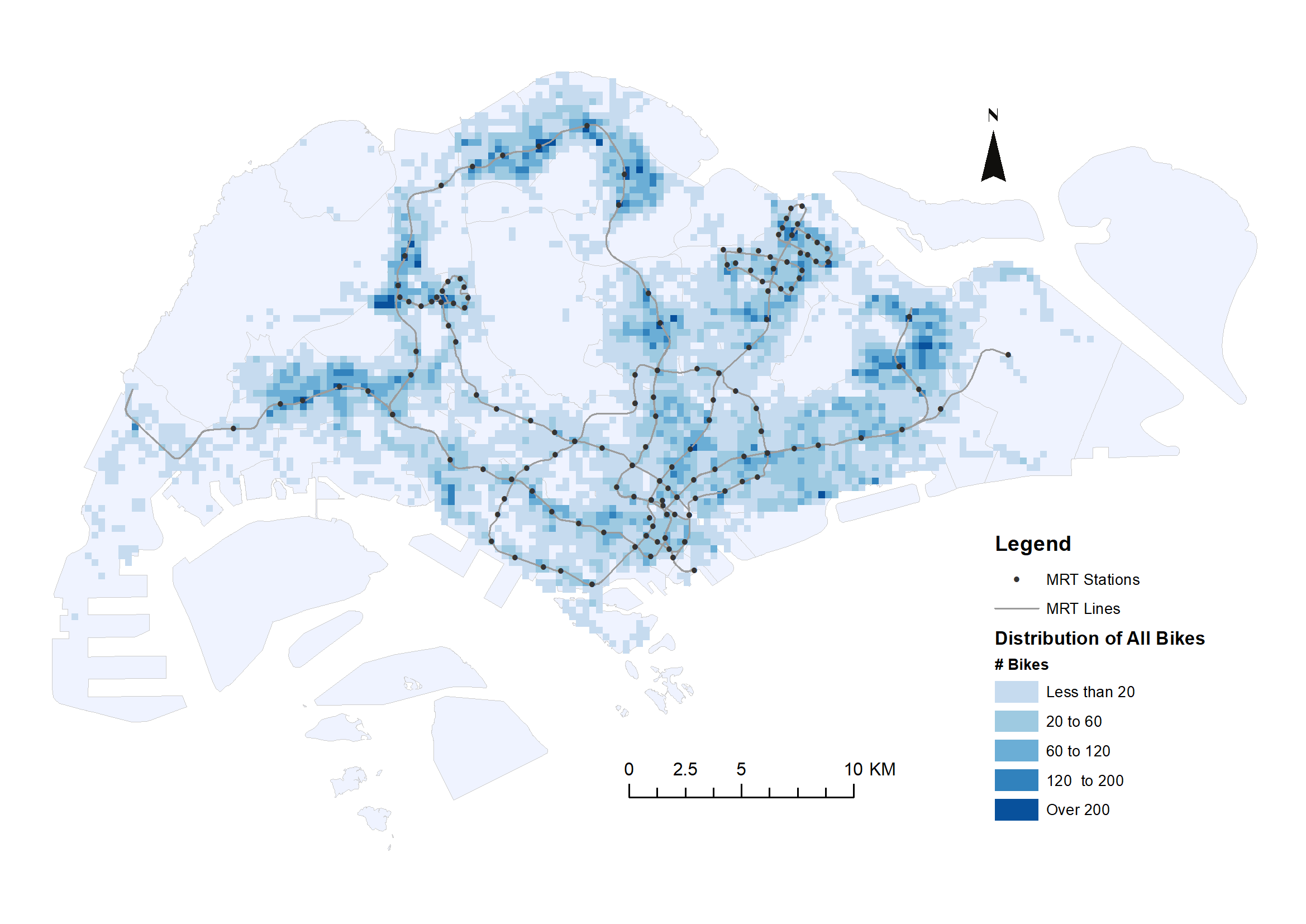}
\caption{Spatial Distribution of All Bikes}
\label{spatial_aday}
\end{figure}

The dockless bikesharing system was first introduced in Singapore in 2017, marking a significant advancement in urban transportation solutions. A distinguishing feature of this system was the integration of GPS sensors on all bicycles, allowing for the continuous tracking of their trajectories. This feature is instrumental in analyzing bike distribution and usage patterns. In our study, we leveraged the anonymous dockless bike GPS data from a bikesharing service provider in October 2017. Figure \ref{spatial_aday} presents an overview of the daily spatial distribution of bikes across Singapore. To provide context, the locations of metro stations are also indicated on the map. As can be seen in the figure, the shared bikes were mainly concentrated around metro stations with a high population density. This strategic placement aligns with the intuitive understanding that areas with larger populations generally yield higher demand, consequently driving increased profitability for operators.

\begin{figure}[ht]
\centering
\includegraphics[width= 0.9 \textwidth]{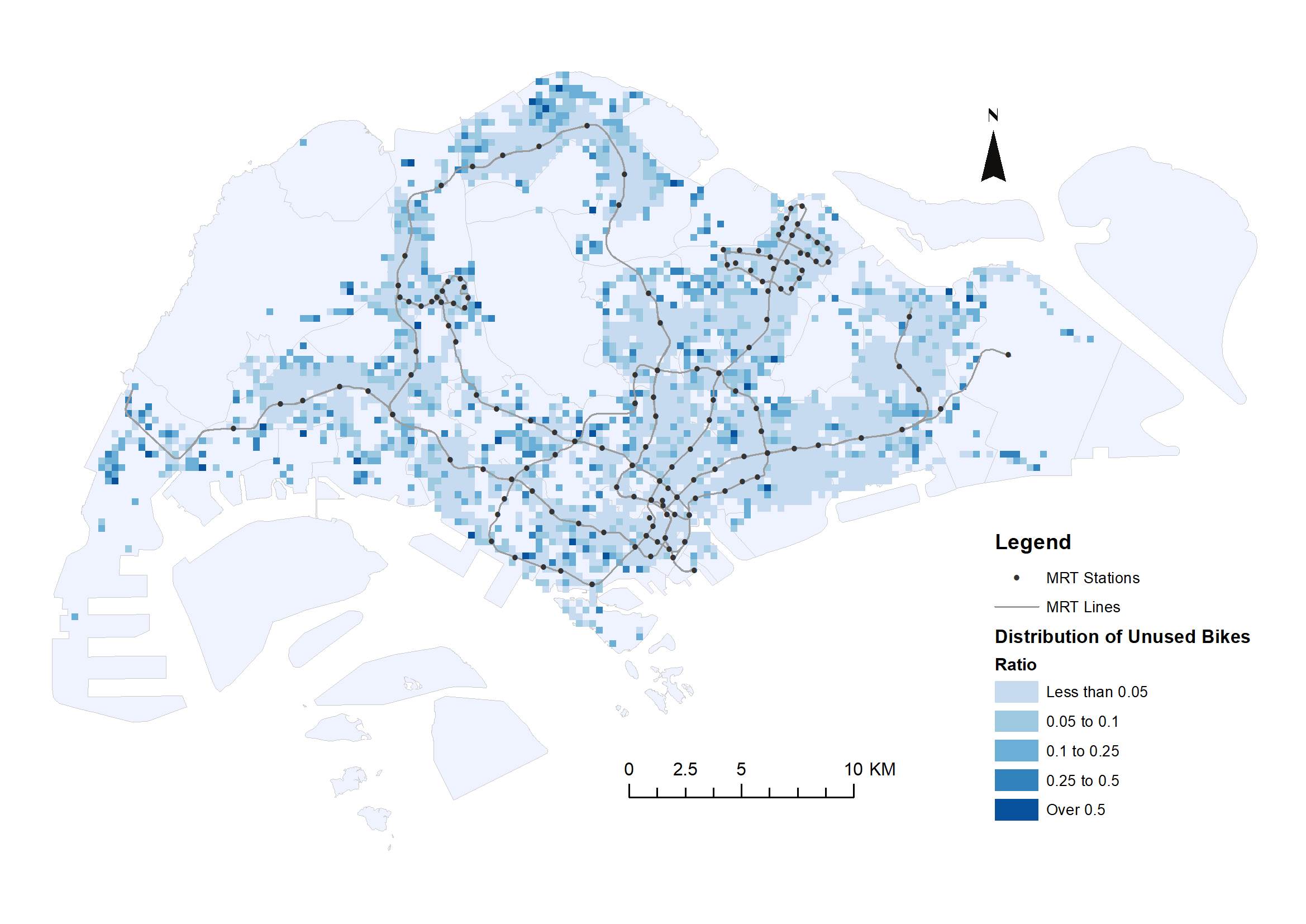}
\caption{Spatial Distribution of the Ratio of Unused Bikes}
\label{spatial_unused}
\end{figure}

Figure \ref{spatial_unused} offers insight into the distribution of the ratio of bikes that remained unused for over one week in October. This ratio serves as an indicator of the oversupply level. Notably, areas with a higher unused bike ratio tended to coincide with areas with a larger number of bikes. Moreover, areas with a high oversupply level were predominantly situated further away from metro stations. This observation underscores the underutilization of bikes in such areas and the necessity for bike fleet relocation. It is imperative to redistribute bikes from remote regions to more central areas to enhance performance within the system.

Bike fleet relocation has been widely implemented in bikesharing systems as an efficient means of fleet management to reach a balance between a limited number of bikes. In previous studies, several models discussing bike fleet relocation strategies have been proposed, including both dock-based and dockless systems. From the modeling perspective, the main difference between a dock-based and dockless bikesharing system lies in the capacity of the regions of interest. While dock-based bikesharing systems feature in a fixed storage capacity, the capacity is usually not a concern when dealing with dockless system.

In this paper, we focus on developing a flow-based optimization model for solving the dockless shared bike relocation problem. The remainder of the paper is organized as follows: In Section \ref{re_work}, we did a short review of related work about bikesharing relocation problem and specified our contribution. In Section \ref{pdn}, we introduced a network flow model for the combined bike relocation and truck trike assignment problem for dockless bikesharing system. We also showed the way to apply the proposed model to a dock-based system by adding capacity constraint. On the basis of the actual demand and distribution of bikes obtained from real-world data, in Section \ref{model_res}, we tested our model in a series of scenarios with various sizes of bike fleets and numbers of truck trikes. The modeling results were further analyzed by performing a series of sensitivity analysis. We also compare the performance of the dock-based and dockless sharing bike systems by solving the proposed model under the two different scenarios. Finally, in Section \ref{conclusion}, we offered some policy suggestions for bikesharing operators and governments for better cooperation to provide a satisfactory public service with reasonable revenues.

\section{Related Works: Modeling Bike Fleet Relocation}\label{re_work}
Numerous studies have delved into the intricate dynamics of bike fleet distribution and relocation strategies in free-flating bikesharing systems. Broadly, relocation problems are classified into two categories distinguished by operating hours:
\begin{enumerate}[noitemsep, topsep=0.5pt, label=(\arabic*)]
    \item Static shared bike relocation at night (neglecting user demand during the relocation period)
    \item Dynamic shared bike relocation (continuous relocation during the daytime while considering the user demand)
\end{enumerate}

Without considering the user demand dynamically during the relocation period, static bicycle relocation problems (SBRP) are easier to formulate and solve than dynamic counterparts due to their lower sensitivity to the user demand. A prevalent approach in solving static bikesharing rebanalcing is to model it under vehicle routing problem (VRP) paradigm. Since VRP is an NP-hard problem, different researchers seek to design effective and efficient algorithms to optimize or approximate solutions. \citet{chemla2013bike} is one of the first to propose modeling the static relocation problem as a VRP problem. The authors were motivated by single vehicle one-commodity capacitated pickup and delivery problem and present a branch-and-cut algorithm for solving its relaxation. Based on this work, \citet{erdougan2015exact} proposed a branch-and-cut algorithm extended by combinatorial benders cuts to get an exact solution for BRP. The algorithm was tested on two instances from literature and demonstrated to be effective. Among other works that seek to solve the same issue, \citet{pal2017free} developed a mixed-integer linear programming (MILP) model for the planning and management of static free-floating bike sharing, using a hybrid nested large neighborhood search with the variable neighborhood descent algorithm. \citet{ho2014solving} presented an n-iterated Tabu search heuristic to solve the static repositioning problem and demonstrate its ability to generate high-quality solutions within limited computing time. While all the aforementioned work sought to improve the computational efficiency from an algorithmic perspective, some work has explored techniques for problem size reduction. \citet{forma20153} designed a three-step heuristics approach for solving the SBRP. Unlike the preivous research that focuses on the relocation optimization based on the initial region discretization, this algorithm starts by the clustering of stations. This step aims to reduce the network size. Then, the problem is solved in a decomposing and solving manner which allows for the application in larger bikesharing network.

While static relocation is easier and more implementable compared to dynamic ones, for a great part of bikesharing systems worldwide, there are usually peak hours during the day when the bike demand is uneven between different stations and distribution is not equal. Dynamic relocation is essential when dealing with uneven demand patterns throughout the day. Under dynamic relocation setting, VRP is still a major paradigm for problem modeling. Unlike static relocation, the target distribution in dynamic relocation needs to be dynamically derived from the demand pattern. The problem is often solved in a predict-and-solve manner, and this process is implemented by running simulation. \citet{schuijbroek2017inventory} utilized the Markov chain to model the state transformation process of a bikesharing system, aiming to reproduce the users’ usage pattern. \citet{caggiani2013dynamic} proposed a bikesharing system simulator to model destination choices and used simulated origin-destination (OD) pairs as the input for optimization models. The demand was predicted by an neural network with a lower and upper bound, These bounds were the main constraints of the optimization model and allowed taking the uncertainty and dynamic features of this system into account. Similar idea can be found in a more recent work by \cite{CAGGIANI2018159}. A clustering step was added before forecasting and optimization as a preliminary step for reducing the problem size and improving computational efficiency. \citet{caggiani2017dynamic} also talked about a dynamic clustering method based on the region activity pattern for the relocation problem in free-floating bikesharing system. As the performance of the prediction model in the simulator has direct influence on the performance of the optimized operation, further specialized research focusing on demand prediction also arises as a result and has been extensively studied. For the sake of more elastic and realistic prediction of bike distribution and demand in bikesharing systems, several researchers have performed comprehensive studies on the relationship between demand and environmental and social factors (\cite{frade2014bicycle, tran2015modeling, faghih2014land, shen2018mobility}). 

More recently, with the development of learning-based control theory, there has been an emerging application of deep reinforcement learning (RL) to bike relocation problem. Reinforcement learning aims to learn to decide optimal relocation operation, which is generally denoted as action, based on the current state of the system, including the demand and vehicle distribution. The quality of the operation is quantified by a manual designed reward function which is related to the achieved following state. The RL agent tries to learn the relationship among the state, action, and reward pair in order to achieve better and better reward over training process. \citet{li2021towards} applied this model to solve the bikesharing relocation problem. The experiment results showed that the framework is able to improve the system's performance. \citet{li2022dynamic} modeled the bikesharing relocation as a muti-agent reinforcement learning problem in which each rebalancing vehicle was treated as an agent. This model gave stable results under various operating condition and showed its robustness.

Apart from using vehicles to redistribute bikes in the system, some researchers have proposed user incentives as an effective means to stimulate users returning bikes to a nonsaturated station \citep{fricker2016incentives,pfrommer2014dynamic}, which is called “user-based” operation. These researchers analyze user decisions in detail and give advice on how to set incentives between stations to encourage user to change their initial destination to a desired one for the operators. While the line of work on incentives mainly focuses on individual bikes, \citet{ghosh2017incentivizing} provides a global view of incentivizing possible trailers, involving multiple bikes and stations, to help meet the repositioning demand. \citet{pan2019deep} trained a novel deep reinforcement learning framework to solve rebalancing problems in dockless bikesharing systems. The goal of this model was to develop an optimal pricing algorithm to incentivize users in congested regions to pick up bikes elsewhere to maximize the service level. This model was extended by \citet{duan2019optimizing}, who considered the possibility of users returning their bikes to alternative destinations. A similar incentive strategy has also been studied by \citet{ji2020does} using a well-defined model. Since user-based model is out of the scope of this paper, interested readers can refer to the above literatures for a more comprehensive overview.

In this paper, our primary focus is on addressing the dynamic relocation problem from an optimization-based perspective. The contributions of this paper are three-fold:
\begin{itemize}
    \item We propose to model the dynamic bike relocation problem by a network flow model. This model is more straightforward to follow and apply compared to the simulation-based approach.
    \item  Based on the proposed model, we conduct a comprehensive sensitivity analysis of system performance concerning factors such as price, cost, and fleet scale. These analyses provide insights for the practical management of the bikesharing system.
    \item We compare the performance of dock-based and dockless bikesharing system by adding the station capacity constraint to the proposed network flow model. The result shows the condition of the transformation from a docked-based to a dockless system from a performance point of view.
\end{itemize}

\section{Problem Statement and Formulation}\label{pdn}

For the purpose of a more clear discussion, we first present the notation that will used to describe the modeling settings and construct the model. Table \ref{tab:notation} shows the complete notation list.

\renewcommand{\arraystretch}{1.5}
\begin{table}[h!]
\begin{center}
    \caption{Notations}
    \begin{tabular}{>{\centering\arraybackslash}m{\textwidth/6}|>{\centering\arraybackslash} m{\textwidth/8}|>{\centering\arraybackslash}m{\textwidth/2}}
    \hline
    \textbf{Type} &\textbf{Notation} &\textbf{Explanation} \\
    \hline
    \multirow{3}{*}{Sets} &$N$ &Set of parking sites for truck trikes\\
    &$A$ &Set of links (representing the shortest paths connecting different parking sites)\\
    &$T$ &Set of time periods, $t\in \text{\{0,1…,T\}}$ \\
    \hline
    \multirow{7}{*}{Parameters} &$c_{ij}^v$ &Cost of transporting a bike from site $i$ to site $j$ by a trike, including loading the bike at $i$ or unloading at $j$\\
    &$c_{ij}^u$ &Cost of operating a trike from site $i$ to site $j$ \\
    &$e_{ij}^d$ &Revenue charged from a rider using a shares bike from site $i$ to site $j$\\
    &$a$ &The capacity of a truck trike\\
    &$n$ &Number of trikes in the bikesharing system\\
    &$r_{ij}^t$  &Number of riders arriving at site $i$ at time $t$ and intending to ride a bike from site $i$ to site $j$\\
    &$P$ &fixed cost for operation per trike normalized to each time unit, i.e., cost of buying trikes divided by their service life\\
    \hline
    \multirow{3}{*}{\shortstack{Continuous\\variables}} &$b_{ik}^{t,t_k}$ &Number of 
    bikes being transported from site $i$ to site $k$ in time period $[t,t_k]$\\
    &$b_i^t$ &Number of bikes being parked at site $i$ at time $t$\\
    &$d_{ij}^{t,t_j}$ &Number of bikes being taken by riders from site $i$ to site $j$ in time period $[t,t_j]$\\
    \hline
    \multirow{2}{*}{Discrete variables} &$w_{ik}^{t,t_k}$ &Number of trikes moving from site $i$ to site $k$ in time period $[t,t_k]$\\
    &$u_i^t$ &Number of trikes waiting at site $i$ at time $t$\\
    \hline
    \end{tabular}
    \label{tab:notation}
\end{center}
\end{table}

In this study, we propose to use a spatiotemporal network to model the movements and stops of bikes and truck trikes in a defined dynamic bikesharing system and analyze their impact on the system performance. In view of the very limited number of truck trikes and the overwhelmingly large number of bikes in the system, our mathematical programming model presented below characterizes trikes as individual identities and describes bikes as continuous flows, which results in a mixed discrete and continuous optimization problem.

Let $n$ be the number of homogeneous trikes reserved in the system, where each trike is associated with capacity $a$, which indicates the maximum number of bikes it can transport in a single movement. Given that there are no docking stations in dockless bikesharing systems, we divided the service area into a number of zones representing naturally formed bike parking sites by customers, where all of these sites typically have a very large limit on the number of bikes parked there. Let $N$ be the set of sites in the selected network. The time horizon is discretized into a number of periods with the same duration. These time periods are indexed from $1$ to $T$ in terms of the temporal order. The number of bikes parked at site $i$ in period $t$ is denoted as $b_i^t$, whereas the number of trikes waiting at site $i$ in period $t$ is $u_i^t$. We then consider a set of decision variables related to moving bikes and trikes between different time periods:
\begin{enumerate}[noitemsep, label=\arabic*)]
    \item The number of bikes being transported from site $i$ to $k$ in time period $[t,t_k]$ is $b_{ik}^{t,t_k}$.
    \item The number of bikes being taken by riders from site $i$ to $k$ in time period $[t,t_j]$ is $d_{ij}^{t,t_j}$.
    \item The number of trikes transporting bikes from site $i$ to $k$ in time period $[t,t_k]$ is $w_{ik}^{t,t_k}$.
    \item The number of trikes moving without any bike from site $i$ to $k$ in time period $[t,t_k]$ is $u_{ik}^{t,t_k}$. 
\end{enumerate}

The gap between $t$ and $t_k$, and that between $t$ and $t_j$ are determined by the distance between the origin-destination (OD) sites and the speed of the bikes or trikes. We assume that the travel time does not depend on the actions performed by a trike at a site. However, this could be a potential limitation of the model that can be partially addressed by considering the average service time.

\begin{figure}[ht]
\centering
\includegraphics[width=0.8 \textwidth]{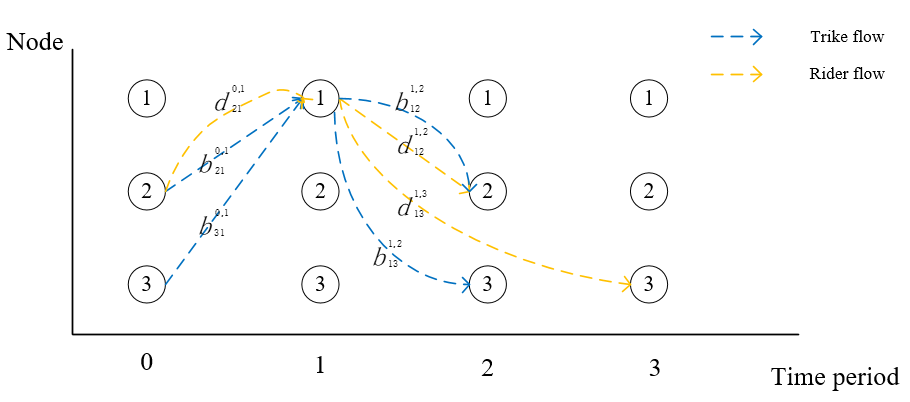}
\caption{Example of a Small Spatiotemporal Network \\
         (Flow for Site 1 at Time 1)
        }
\label{network_example}
\end{figure}

A toy example with three sites and three time periods is illustrated by a time–space graph shown in Figure \ref{network_example}. Only the flows of bikes and trikes that will affect the number of bikes at Site 1 at Time 1 are shown. The flow-in and flow-out numbers of bikes at Site 1 consist of two parts: trike transportation and rider usage. The relationship between the number of bikes at two consecutive moments can be precisely formulated using these flows. As in this example, trike and rider flows from Site 1 to Site 3 vary because of the different speeds of the bikes and trikes. Moreover, when deciding the traveling time between two sites, we need to focus only on the distance between the two sites, not the routes. For example, even if Site 2 is located in the route from Site 1 to Site 3, this does not influence the results. Note that if the distance from Site 1 to Site 2 in the real network takes half a time period to cover by bike, it should be rounded up to 1 in the spatiotemporal network.

\subsection{Problem formulation}\label{prob_form} 
In this section, we propose a MILP model to characterize the combined bike relocation and truck assignment problem. The goal of this model is to find an optimal dynamic operations scheme for assigning truck trikes to relocate bikes between parking sites to maximize the system operations profit for a near-future short-term period. The operations profit defined here is simply the revenue charged to riders to use the bikes minus the bike relocation cost (including the loading, transportation, and unloading costs) and trike operations cost over the analysis period. Other indirect costs, such as vehicle depreciation cost and trike drivers’ salaries and welfare, are better evaluated on an intermediate-term or long-term scale and, hence, should not be taken into account by the operations model. Of course, the model focuses on the operations of relocating bikes by a limited number of trikes to redistribute the bike supply to serve the rider demand over different sites to the maximum extent. The major operational decisions of this model are, during each time interval, how many bikes should be relocated from one site to another and how to assign trucks to transport those bikes between different sites.

To capture the most essential features of the described bikesharing operations problem while avoiding unnecessary modeling difficulties, the proposed model implies the following individual behavior and system modeling settings or assumptions:
\begin{enumerate}[noitemsep, label=\arabic*)]
    \item The entire system analysis period is discretized into a number of time intervals with an equal duration.
    \item The number of bike usage requests for each site pair during each time interval can be accurately predicted.
    \item All riders can make their bike-riding trips from their pickup sites to dropoff sites directly without detouring to other places or conducting other activities over the course of their trips.
    \item While bikes and trikes move at different speeds in the network and their movement speeds are subject to congestion conditions in the network, it is assumed that their travel times between sites can be deterministically and reliably predicted.
    \item The charge for a rider's bike usage is a function of his or her travel time.
    \item The movements of bikes and riders in the network form aggregate continuous flows, whereas the movement of trikes us counted individually.
    \item Riders leave their pickup site if they find no bike available at the time they arrive at the site; in other words, the riders are not willing to wait for incoming bikes.
    \item Both the loading and unloading times for transferring bikes between a parking area and a truck are ignored, which, if necessary, could be considered without adding too much complexity to the current model.
    \item The parking area of any site for bikes is sufficiently large so that no capacity constraint needs to be imposed on any site.
\end{enumerate}

By incorporating all assumptions listed above, a neat MILP model embedded with a spatiotemporal network structure with multiple vehicles/commodities (including riders, bikes, and truck trikes) reads as follows: 

\begin{maxi!}
    {}{\sum_{t}\sum_{i}\sum_{j}e_{ij}^{d}d_{ij}^{t,t_j}-\sum_{t}\sum_{i}\sum_{k}c_{ik}^{u}w_{ik}^{t,t_k}-\sum_{t}\sum_{i}\sum_{k}c_{ik}^{v}b_{ik}^{t,t_k}-nP \label{eq:Obj}}
    {\label{eq:Eq1}}{}
    \addConstraint{b_{i}^{t+1}}{=b_{i}^{t}-\sum_{k}b_{ik}^{t,t_k}-\sum_{j}d_{ij}^{t,t_j}+\sum_{k}b_{ki}^{t,t_k}+\sum_{j}d_{ji}^{t_j,t}, \; \forall i\in N, t \in T \label{eq:Con1}}
    \addConstraint{u_{i}^{t+1}}{=u_{i}^{t}-\sum_{k}w_{ik}^{t,t_k}+\sum_{k}w_{ki}^{t_k,t}, \; \forall i\in N, t \in T \label{eq:Con2} }
    \addConstraint{b_{ik}^{t,t_k}}{\leq aw_{ik}^{t,t_k}, \; \forall i,k \in A, t \in T \label{eq:Con3} }
    \addConstraint{d_{ij}^{t,t_j}}{\leq r_{ij}^{t,t_j}, \; \forall i,j \in A, t \in T \label{eq:Con4} }
    \addConstraint{b_{i}^{t}}{\geq 0, \; \forall i\in N, t \in T \label{eq:Con5}}
    \addConstraint{u_{i}^{t}}{\in \{ \, 0,1,2,\ldots \}\, , \; \forall i\in N, t \in T \label{eq:Con6}}
    \addConstraint{b_{ij}^{t,t_j}}{\geq 0, \; \forall i,j\in A, t \in T \label{eq:Con7}}
    \addConstraint{d_{ij}^{t,t_j}}{\geq 0, \; \forall i,j\in A, t \in T \label{eq:Con8}}
    \addConstraint{w_{ik}^{t,t_k}}{\in \{ \, 0,1,2,\ldots \}\, , \; \forall i,k\in A, t \in T \label{eq:Con9}}
\end{maxi!}
    
The objective function (\ref{eq:Obj}) maximizes the system profit,  which is calculated as system revenue charged to riders minus the bike relocation cost and truck operations cost. The constraints (\ref{eq:Con1}) and (\ref{eq:Con2}) bike flow conservation and trike number conservation at parking sites, respectively. The constraint (\ref{eq:Con3}) guarantees that the total number of bikes transported by trikes does not exceed the total capacity of trikes. The constraint (\ref{eq:Con4}) sets the demand satisfaction bound (i.e., the number of bikes used by riders or the number of riders picking up bikes is less than the number of riders arriving at any site over any time interval). The constraints (\ref{eq:Con5}) and (\ref{eq:Con6}) ensure that all bike flow variables are nonnegative. The constraints (\ref{eq:Con7}), (\ref{eq:Con8}), and (\ref{eq:Con9}) ensure that all bike flow variables are nonnegative and that all trike numbers are nonnegative integers subject to the total number of trikes in the network.

\subsection{From dockless to dock-based system}
Dock-based bikesharing systems have been around for much longer than dockless systems and are considered a typical scenario of issues related to bikesharing systems. A large number of studies have focused on the rebalancing problem of dock-based bikesharing systems. Most of these studies have regarded this problem as a vehicle routing problem and highlighted the importance of proper routes for every vehicle in the system to fill the gap between the current system state and expected state \citep{chemla2013bike,liu2016rebalancing,pfrommer2014dynamic}. Therefore, our aim was to use our model to solve and understand this problem from a different perspective. In our model, because we have discretized space in dockless bikesharing systems, we only need to add capacity constraints for each site to solve the rebalancing problem for dock-based bikesharing systems. With regard to our model, we concentrated on the revenue of the whole system, although previous studies have mainly focused on improving the customer service, that is, improving the probability of riders finding a bike to ride and having an empty dock to return the bike to whenever and wherever they want.

We then add the following constraint to the model proposed in Section \ref{prob_form}:
\begin{equation}
    {b_{i}^{t}}{\leq cap, \;  \forall i\in N, t \in T}
\end{equation}
where cap denotes the given capacity for each site. With this capacity constraint, we can use the model to solve relocation problem under docked-based system assumption. We will use this model to analyze the transition condition between the two type of systems later in the case study.
\section{Case Study and Analysis}\label{model_res}
In this section, we apply the above-proposed model to the bike relocation problem by utilizing the real operational data in Singapore. After a brief introduction of the dataset, we construct, solve, and analyze the result of a base case to show the effectiveness and practicality of the model. Subsequently, we conduct a sensitivity analysis of the system performance to the system configuration factors such as the price, cost, and fleet size. Additionally, we explore the implications of different operational assumptions by comparing outcomes between dock-based and dockless systems from the proposed model. The objective is to derive valuable insights regarding the system features and operations, which is further translated into practical policy suggestions.
\subsection{Data preparation and experiment setup}
Singapore is divided into several districts, with a district called Yishun located in the north. Being a densely populated residential area, Yishun has a high demand for bikesharing services, which makes bike relocation necessary. As a result of such a high usage frequency of shared bikes, Yishun has abundant historical data of bikesharing systems for conducting numerical experiments.

\begin{figure}[H]
\centering
\includegraphics[width= 0.9 \textwidth]{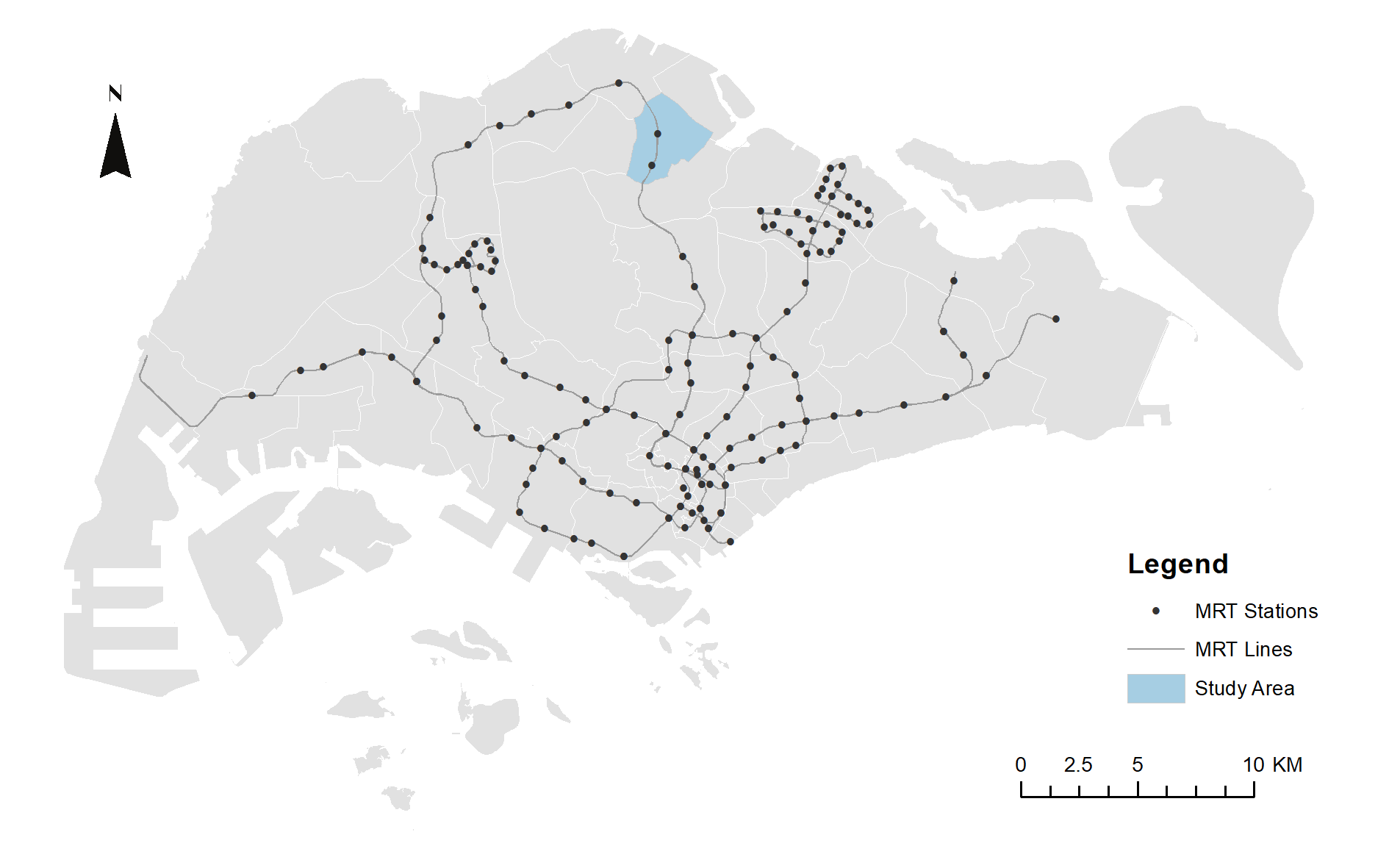}
\caption{Geographical Location of the Study Area}
\label{yishun}
\end{figure}

In a study by \citet{shen2018mobility}, the authors demonstrated that most of the bike rides occur within each district locally in Singapore, showing a geographically clustered usage pattern. Such a small case study performed within the area of Yishun can provide us with enough insights into the operation details and influence of different factors on the system operations. Therefore, we decided to select Yishun as our case study area, which is highlighted in Figure \ref{yishun} to provide a more intuitive idea of its location. In the following numerical experiments, we divided the area of Yishun into equally sized land cells with a size of 300 × 300 m, with each cell standing for a single site in our network model.

In the numerical experiments, we utilized the data of weekdays during the morning peak hours from 7:00 to 10:00. It was assumed that relocation can be performed at any time during these hours. For the sake of easier modeling, the time was discretized into an interval of 5 min, and the bike trips were aggregated to intervals in each cell. Theoretically, these aggregated trips are the achieved demand rather than the true demand because the true demand should also contain the potential unachieved demand. However, since there was an oversupply of shared bikes during the observation period in Yishun, we assume the achieved demand is a fair approximation of the true demand in this case. As dockless bikes sometimes appear in random locations, to control the problem size within a tractable scale, cells with a total demand of less than 10 trips during these three hours were filtered out.

Regarding the price and cost computation mechanism, we assume the profit from riders was computed as a function of the time duration of the trips. The distance between sites was computed according to the road network distance between the sites, the time duration was subtracted from the trip data, and the parameters for this price model were obtained from the actual price scheme. The cost of relocation consists of the bike abrasion, labor cost, and the travel cost of each trike corresponding to gas consumption. The inital distribution was computed from the data and used as the input for solving the model. 

Based on the above processing and definitions, we commenced the experiment by computing a base case with 10 sites to validate the reliability of the model. Then, we scale up the base case and perform a sensitivity analysis of the important aspects of bikesharing relocation problems, including the number of trikes and bikes, operation cost, and price of using bikes. The provide an overview, all the parameters studied in the sensitivity analysis as well as their value sets are summarized in Table \ref{sens_case}. While conducting these sensitivity analyses, we only changed the value of one parameter at one time while keeping the others unchanged as the default value. The default values for the number of trikes, number of bikes, cost, and price were 2, 899, 0.5, and 1, respectively. All the other unmentioned parameters were set to be the same as in the base case.

\begin{table}[ht]
\centering
\caption{Scenarios for the Sensitivity Analysis}
\begin{tabular}{c|c}
\hline
\textbf{Parameters} & \textbf{Value set} \\
\hline
number of trikes       & 0,1,2,3,4,5        \\
number of bikes        & 225,488,602,899     \\
cost (in S\$)            & 0.4,0.6,0.8,1      \\
price (in S\$)           & 1,1.5,2,2.5,3     \\
\hline
\end{tabular}
\label{sens_case}
\end{table}

Next, we delve into more details about the different experiments. All the experiments presented in the following part were solved in Gurobi 9.0.2 \citep{gurobi2018gurobi} on a computer with an Intel Core i7-9750H processor.

\subsection{Base case: Analysis of relocation operations}
The values of some basic parameters that we used for the base case experiment are presented as follows. Without explicitly mentioning, these values are kept constant across different experiment cases.
\begin{enumerate}[noitemsep, label=\arabic*)]
    \item The profit is S\$1 every 15 min or less than 15 min.
    \item The cost of operating a trike is S\$0.4/km.
    \item The cost of loading and unloading a shared bike is S\$0.5.
    \item The capacity of each trike for shared bikes is 20.
    \item The number of trike is 3.
    \item The cost of buying one trike scaled to every hour is S\$70.
\end{enumerate}

\begin{figure}[ht]
    \centering
    \includegraphics[width=0.5 \textwidth]{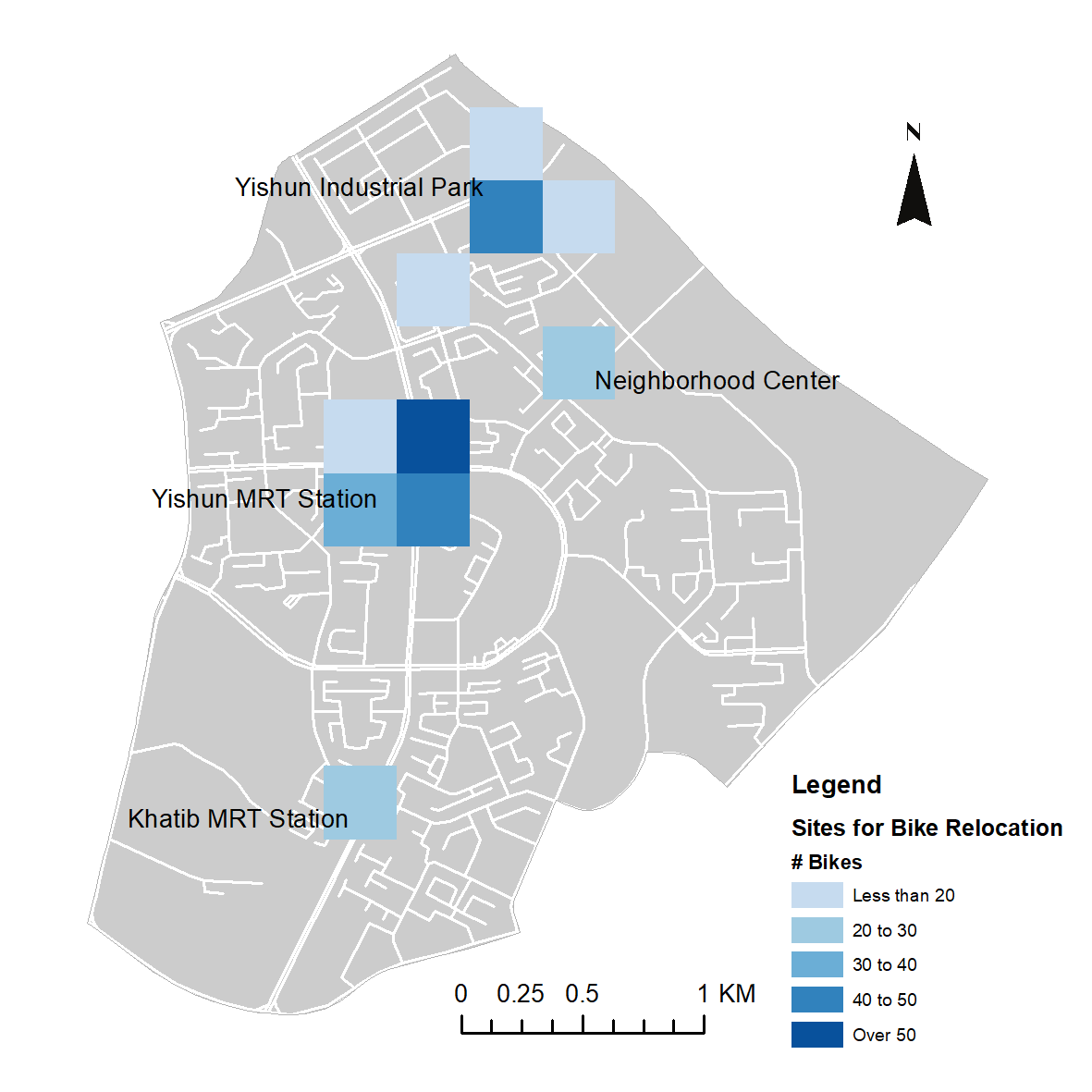}
    \caption{Selection of sites in the Study Area}
    \label{selected10}
\end{figure}

In this section, we use a time-space graph to visualize the exact operation of trikes for relocation given by our model. To prevent the graph from becoming too complex, we intentionally control the scale of the case to be proper for a full visualization, comprising 10 sites and three trikes. Without the loss of generality, those 10 sites are picked to be located around Yishun and Khatib MRT stations, Yishun Industrial Park, and a neighborhood center, which are areas with the highest bike usage demand (see Figure \ref{selected10}). 

\begin{figure}[htb!]
    \centering
    \includegraphics[width=1.0 \textwidth]{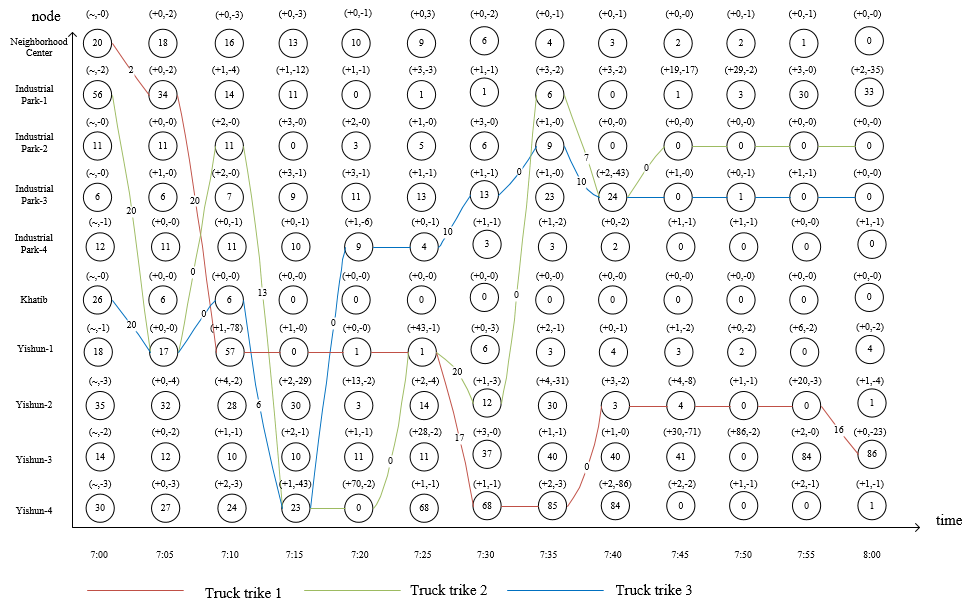}
    \caption{Visualization of Relocation Operations}
    \label{visualization}
\end{figure}

\begin{figure}[htb!]
    \centering
    \includegraphics[width=0.5 \textwidth]{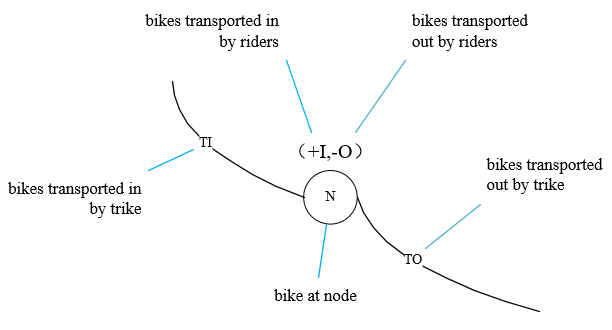}
    \caption{Elements Illustration}
    \label{element_illus}
\end{figure}

 Figure \ref{visualization} shows the relocation operations during the first 12 intervals. The horizontal axis on the spatiotemporal graph is the time horizon and the vertical axis the 10 site names, among which Yishun-1 to Yishun-4 denote the four selected sites around the Yishun MRT station, and Industrial Park-1 to Industrial Park-4 denote the four selected sites around the Yishun Industrial Park. Note that the proximity of two sites on y-axis does not indicate a closer geographic distance. A more detailed explanation about the implications of the circles and arcs is displaied in Figure \ref{element_illus}. The number in each circle indicates the number of bikes parked at each site at the corresponding timestamp, and the tuple above each circle stands for the bike inflow and outflow. In the following explanation, every circle, which includes attributes of time and site, is called a state. 

In Figure \ref{visualization}, the used bikes shows a pattern of intensive use at a certain state, for example, at Time 7 and Site Yishun-2, and a relative stochastic use at other states. Moreover, since the neighborhood center is within walking distance to the Yishun Industrial Park and Yishun MRT station, people living near the neighborhood center do not need to use shared bikes much, so both the inflow and outflow of bikes are small in this area on the graph. If we take a closer look on the relocation arcs, we can find that relocation mainly occurs between different sites in the four areas: Yishun Industrial Park, neighborhood center, Khatib MRT station, and Yishun MRT station. Notably, relocation occurs most frequently between the Yishun Industrial Park and Yishun MRT station. This result can be explained by people's routine: Yishun MRT station is the nearest station to the park. People who need to go to the park by subway during the peak hours mostly get off at the Yishun MRT station and then walk or ride a bike to the park. To go back, they usually walk or ride a bike from the park to the station. Intensive use rate is expected at the Yishun Industrial Park and Yishun MRT station, which is the reason for the frequent relocation between these two areas. In contrast, the Khatib MRT station is further from the neighborhood center and the park, the two important centers of Yishun, than the Yishun MRT station, so the Khatib MRT station has a lower shared bike use rate and hence requires fewer relocation operations. 

A more detailed analysis of the relocation operation reveals that it occurs mainly between the Yishun MRT station area and Yishun Industrial Park from 7:00 to 7:30 and within the Yishun MRT station area from 7:00 to 8:00. Importantly, the relocation operation pattern observed from 7:00 to 7:30 occurs because the departure demand at the Yishun MRT station from 7:30 to 8:00 is greater than that from 7:00 to 7:30. Under this case, relocation needs to be performed in advance to prepare for the increasing demand during the peak hours. From 7:30 to 8:00, many riders arrive at Yishun-1 and many others need to depart from Yishun-2 and Yishun-4, so relocation occurs within the Yishun MRT station area during this period to meet this demand.

Overall, the trajectory of each trike from the graph is intuitively reasonable without irregular operations, such as back-and-forth relocation between certain site pairs or completely idle time. In conclusion, the relocation results obtained by our model are reasonable and implementable.

\subsection{Sensitivity analysis}
In this section, we present the results of the sensitivity analysis regarding the number of trikes and bikes, operation cost, and price on revenue. We first analyze infrastructure scales, including the numbers of trikes and bikes, in two different scenarios to provide insights into the proper number of trikes and bikes in different situations. Then, we perform a sensitivity analysis of two core factors, cost and price. These two factors are directly related to revenue and thus need to be deliberately evaluated to ensure the system's profitability. Assessing the influence of these four factors on revenue is the basis for a successful business. Table \ref{studied scenarios} summarizes the values of these parameters in each scenario.

\begin{table}[ht]
\centering
\caption{Studied scenarios}
\label{studied scenarios}
\begin{tabular}{c|cccc}
\hline
\multicolumn{1}{l|}{}                                                                      & \textbf{number of trikes} & \textbf{number of bikes} & \textbf{cost in S\$} & \textbf{price in S\$} \\ \hline
\multirow{6}{*}{\begin{tabular}[c]{@{}c@{}}Sensitivity\\ to\\ number of\\ trikes\end{tabular}} & 0                      & 899                   & 0.5               & 1                  \\
                                                                                           & 1                      & 899                   & 0.5               & 1                  \\
                                                                                           & 2                      & 899                   & 0.5               & 1                  \\
                                                                                           & 3                      & 899                   & 0.5               & 1                  \\
                                                                                           & 4                      & 899                   & 0.5               & 1                  \\
                                                                                           & 5                      & 899                   & 0.5               & 1                  \\ \hline
\multirow{4}{*}{\begin{tabular}[c]{@{}c@{}}Sensitivity\\ to\\ number of\\ bikes\end{tabular}}  & 2                      & 225                   & 0.5               & 1                  \\
                                                                                           & 2                      & 488                   & 0.5               & 1                  \\
                                                                                           & 2                      & 602                   & 0.5               & 1                  \\
                                                                                           & 2                      & 899                   & 0.5               & 1                  \\ \hline
\multirow{4}{*}{\begin{tabular}[c]{@{}c@{}}Sensitivity\\ to\\ cost\end{tabular}}           & 2                      & 899                   & 0.4               & 1                  \\
                                                                                           & 2                      & 899                   & 0.6               & 1                  \\
                                                                                           & 2                      & 899                   & 0.8               & 1                  \\
                                                                                           & 2                      & 899                   & 1                 & 1                  \\ \hline
\multirow{5}{*}{\begin{tabular}[c]{@{}c@{}}Sensitivity\\ to\\ price\end{tabular}}          & 2                      & 899                   & 0.5               & 1                  \\
                                                                                           & 2                      & 899                   & 0.5               & 1.5                \\
                                                                                           & 2                      & 899                   & 0.5               & 2                  \\
                                                                                           & 2                      & 899                   & 0.5               & 2.5                \\
                                                                                           & 2                      & 899                   & 0.5               & 3                  \\ \hline
\end{tabular}
\end{table}

While performing our sensitivity analysis of the optimal model solutions, we also performed a sensitivity analysis of the consumer satisfaction level. While revenue ensures short-term profitability, the satisfaction level influences the sustainability of bikesharing systems in the long run and is directly related to social welfare. The consumer satisfaction level is defined as the ratio of the total number of satisfactory trips to the total demand of dockless shared bikes. Here, the number of satisfactory trips refers to the number of riders who find bikes to ride whenever they want.

\begin{figure}[htb!]
    \centering
    \includegraphics[width= 0.5\textwidth]{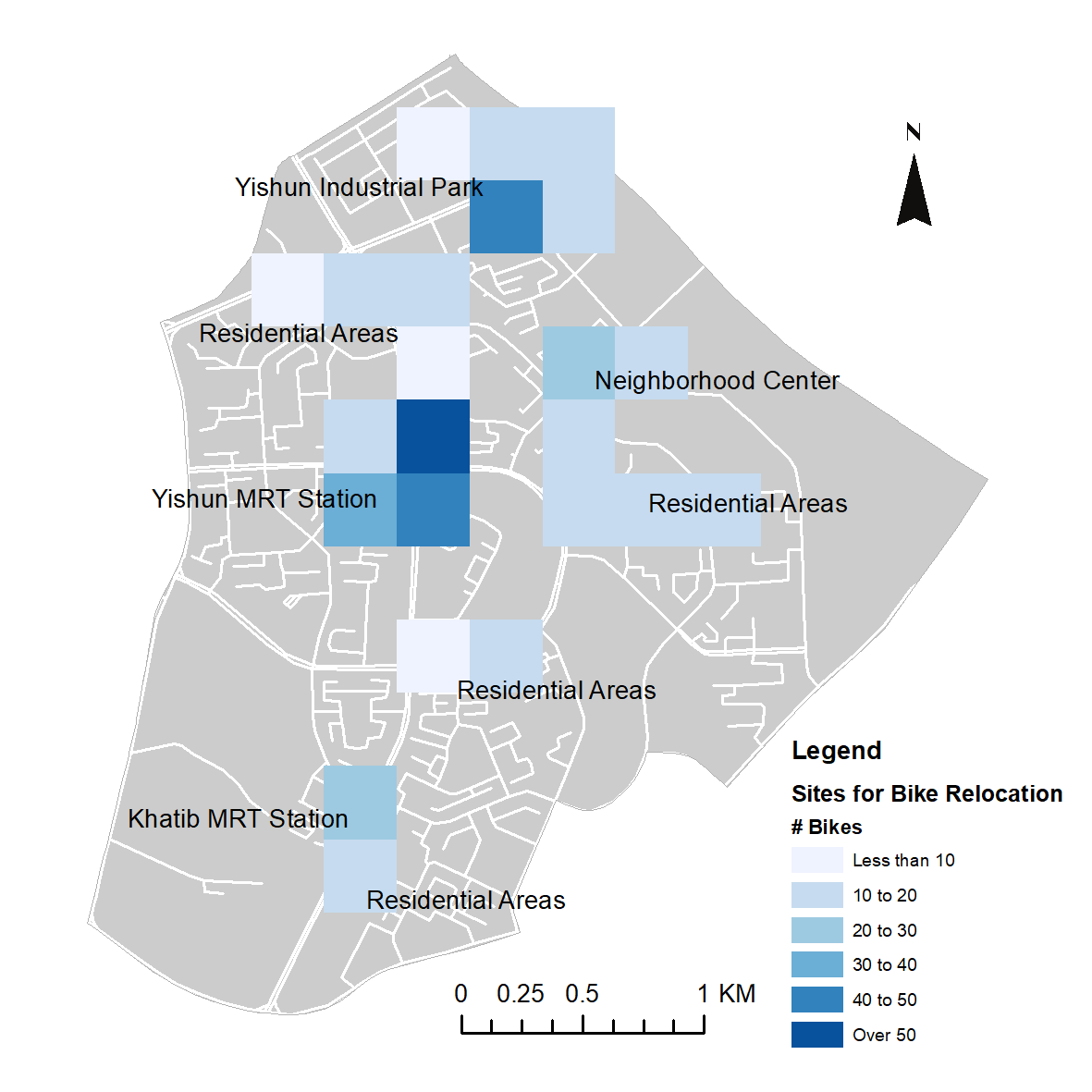}
    \caption{An Extended Selection of sites in the Study Area}
    \label{selected}
\end{figure}

As an extension of the base case, in the sensitivity analysis, we selected the sites with the highest numbers of trips in Yishun. Refer to Figure \ref{selected} for the locations of these sites of interest.

\subsubsection{Sensitivity to infrastructure scale}
Generally, in a bikesharing system, there are two types of decisions that the operators need to carefully make:
\begin{enumerate}[noitemsep, label=\arabic*)]
    \item The bikesharing operators seek to use as fewer trikes as possible for bike relocation with the current bikes to maximize the profit.
    \item The bikesharing operators seek to decide the least number of bikes they need to put into the market with a limited number of trikes to maximize their profit.
\end{enumerate}

\begin{figure}[ht]
\centering
\subfloat[taking fixed cost out of consideration in revenue]{\includegraphics[width=0.5\textwidth]{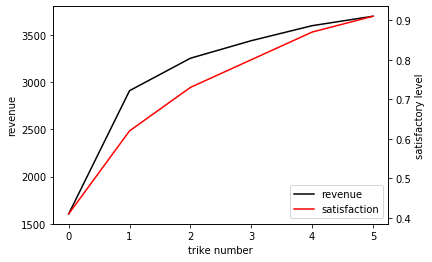}\label{sens_trike_no}}
\hfil
\subfloat[taking fixed cost into consideration in revenue]{\includegraphics[width=0.5\textwidth]{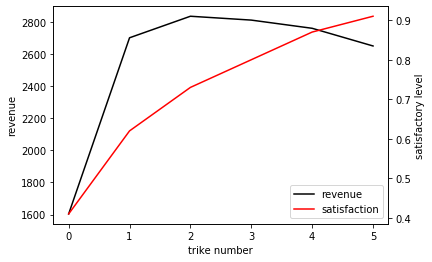}\label{sense_trike_with}}
\hfil
\caption{Sensitivity to Number of Trikes}
\label{sens_trike}
\end{figure}

To resolve the first scenario, we analyzed the profit sensitivity to the number of trikes with a fixed number of bikes. Figure \ref{sens_trike} shows the results and all cases were calculated with a computation time limit of 30 min. As displaied in Figure \ref{sens_trike_no}, when the cost of purchasing trikes is ignored, the revenue continues monotonically with the number of trikes with a decreasing increase speed, and leads to a converging revenue in the end. In Figure \ref{sense_trike_with}, the revenue and satisfaction level increase when the number of trikes increases. They reach the peak when the number of trikes is 2, and fater that, the revenue starts to decrease. The peak observed in Figure \ref{sense_trike_with} exists because relocation can make up for bike insufficiency at each site, thus enabling to meet a higher demand; hence, the presence of more trikes indicates a higher bike usage rate in the system. However, as buying more trikes indicates a higher fixed cost for relocation, it is not rational to purchase many trikes in the system to obtain a higher revenue. A proper number of trikes should be deliberately decided to reach the peak revenue.

\begin{figure}[ht]
    \centering
    \includegraphics[width=0.5 \textwidth]{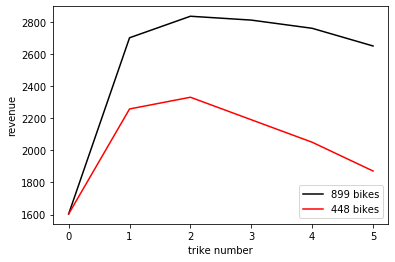}
    \caption{Compare Between Different Numbers of Bikes}
    \label{trike_compare}
\end{figure}

To investigate the influence of the number of trikes on the optimal revenue with different bike fleet sizes, we conducted another experiment with only half the number of bikes as in the real case. Figure \ref{trike_compare}reveals that the revenue generated by half bike fleet size with two trikes case was higher than that of with 899 bikes and no trikes. This observation underscores a pivotal conclusion: by choosing the proper number of trikes for relocation, the profit obtained with fewer bikes in the system can approach or even exceed that obtained with a greater number of bikes. This result validates that it is not essential to amass an extensive fleet of bikes to maximize profit. Instead, strategically implemented relocation can effectively serve as a substitution to meet trip demands and consequently improve revenues.

\begin{figure}[ht]
    \centering
    \includegraphics[width=0.5 \textwidth]{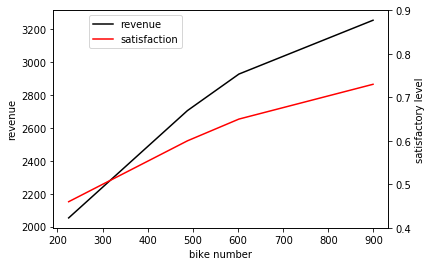}
    \caption{Sensitivity to Bike number}
    \label{sens_bike}
\end{figure}

To address the second scenario, we analyzed the revenue sensitivity to the number of bikes. In this part of the experiments, the number of trikes was first fixed at two while the bike fleet size was systematically varied. The number of bikes is changed by scaling down the initial distribution across site proportionally. Figure \ref{sens_bike} illustrates the revenues obtained under these varying bike fleet sizes. The satisfaction level increased monotonically with the number of bikes. This is intuitive as from the riders’ perspective, a higher number of bikes indicates a better service quality. Although a positive relationship between the number of bikes and revenue existed, it did not imply that a larger bike fleet size was inherently beneficial. To grasp the full context, we must consider the underlying assumptions of this scenario. The second scenario reflects the initial phase of luanching a new bikesharing system. During this stage, the cost of purchasing bikes is a pivotal budgetary factor. However, it is important to note that the bike purchasing cost is not explicitly incorporated into the model's objective in formular \ref{eq:Obj}. Since the cost of this part is fixed since the number of bikes is known for every case, the true final profit can be calculated by directly subtracting the cost of purchasing the bikes from the current revenue. It is worth mentioning that because we only consider the relocation problem in a limited time window, we normalized the fixed cost to each time unit before subtracting. For instance, if we assume that the average life cycle of a shared bike is one year, the cost of purchasing bikes will be evenly distributed over each day, which would amount to S\$1.5. Table \ref{sens_profit} presents the results with the updated revenue. Notably, no monotonic relationship can be observed between the number of bikes and profit. In fact, the final profit is influenced by many factors, including the price of the bike, the life cycle of bikes, and the riders’ treatment pattern for the shared bikes. Given that profit and satisfaction levels do not consistently increase with the number of bikes, determining the optimal number of bikes in the market should be based on a comprehensive analysis process. The result may vary depending on the operating periods, considering the dynamic nature of system features throughout the day. Moreover, the decision may also evolve during different operational phases of the business. For example, during the initial launch of a bikesharing program, the satisfaction levels tends to be crucial as they play a vital role in gaining customers’ trust. However, as the program progresses into a stable operation phase, profitability becomes a more and more important factor to consider.

\begin{table}[ht]
    \centering
    \caption{The Sensitivity of Profit to Number of Bikes}
    \label{sens_profit}
    \begin{tabular}{c|c|c}
    \hline
    \textbf{Number of bikes} &\textbf{Number of trikes} &\textbf{Profit}\\
    \hline
    899 &2 &1905\\
    602 &2 &2024\\
    488 &2 &2033\\
    225 &2 &1715\\
    \hline
    \end{tabular}
\end{table}

\subsubsection{Sensitivity to operation cost}
\begin{figure}[ht]
    \centering
    \includegraphics[width=0.5 \textwidth]{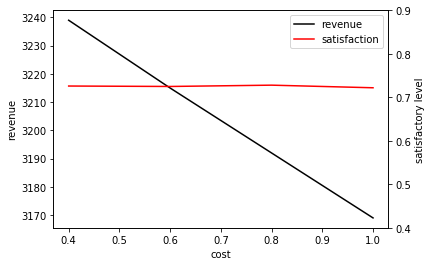}
    \caption{Sensitivity to Cost}
    \label{sens_cost}
\end{figure}

We now evaluated the effects of changes to the relocation cost on the operation results. Each experiemet in this part featured a consistent system configuration with 899 bikes and two trikes. Figure \ref{sens_cost} shows the results of these experiments. It was observed that when the operation cost doubled, the impact on revenue was relatively modest, with a reduction of less than 2\%. Equally significant was the fact that the satisfaction level remained relatively stable despite changes in operational costs. This effect can be attributed to the nature of the cost increases, which primarily pertains to the operation of trikes. Given that each trike serves a fleet of bikes, any increased cost is distributed among this number. Consequently, the influence of the operational cost does not exhibit a pronounced effect in comparison to changes in the scale of infrastructure.

\subsubsection{Sensitivity to usage price}
As the pricing scheme is one of the core steps for a successful business, we also discuss how price affects operations and results. Similar to what we did in previous cases, we changed the price on the basis of the case with 899 bikes and two trikes in the system.

\begin{figure}[ht]
    \centering
    \includegraphics[width=0.5 \textwidth]{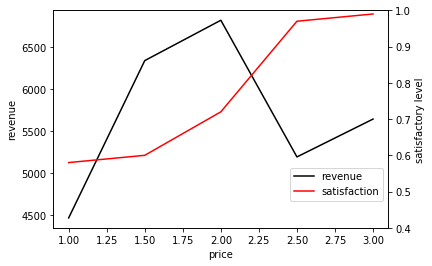}
    \caption{Sensitivity to Price}
    \label{sens_price}
\end{figure}

In the sensitivity analysis with a fixed price, we defined the satisfaction level as the ratio of the met demand of trips to the total demand of trips in each case. We also assumed that the customers are not sensitive to price changes. However, it's imperative to acknowledge that in the real world, customer choices are intricately intertwined with pricing. In practice, as prices fluctuate, the demand for services may vary accordingly, with demand decreasing as prices increase. This inherent relationship between price and demand necessitates careful consideration when modeling the effects of price changes. To address this dynamic, we integrated a well-established price-demand relationship, as presented by \citet{kaviti2019assessing}, into our analysis. Figure \ref{sens_price} provides a visualization of our findings. The result reveals that the revenue trend does not consistently align with price change due to demand fluctuations. However, the satisfaction level displays a notable pattern. As the price increases, the satisfaction level rises. This is because, with the same number of bikes in the system, the demand decreases with the increase in price. This seemingly paradoxical relationship can be elucidated. With a consistent number of bikes in the system, as prices increase, the demand tends to decrease. This decrease in demand results in a higher number of available bikes per customer, effectively increasing the probability that a riding demand will be served. Consequently, the satisfaction level, defined as the exact probability of a riding demand being met, exhibits a positive correlation with price increases.

\subsection{Comparison with dock-based system}

\begin{figure}[ht]
    \centering
    \includegraphics[width=0.5\textwidth]{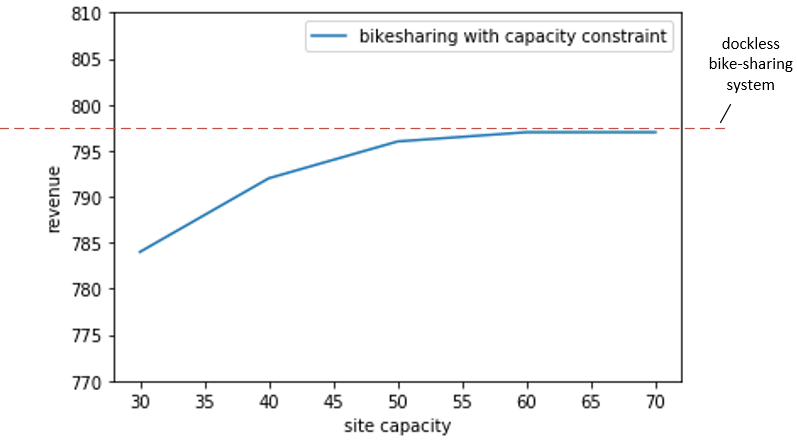}
    \caption{Revenue Changes with Given Capacity}
    \label{revenue_cap}
\end{figure}

In Section \ref{pdn}, we present the model under both dockless and dock-based settings. In this section, we compare their optimization results by setting the initial state and demand at each site the same for the numerical experiments of the two systems. Based on the analysis of the dataset, bikes in dockless bikesharing systems exhibit an apparent uneven distribution pattern. At certain sites, the number of bikes is almost 10 times that in others. We need to set the capacity high enough to make the initial distribution meet the capacity constraint if we want to use the original data. However, setting a high capacity at the beginning will make all the bikes at different sites meet the capacity constraint automatically thereafter, which will nullify the constraint. Hence, to solve this problem, we scaled down the initial number of bikes at each site to less than 30 and increased the capacity from 30 to 70. We selected 10 sites and set three trikes as well as other parameters the same as in the base case in all experimental cases. The results are shown in Figure \ref{revenue_cap}. 

With the increase in the capacity at each site, the revenue of the dock-based bike system increases and becomes equal to that in dockless bikesharing systems after reaching a capacity of 60. This result proves that by constructing a proper number of docks at each site, dock-based bikesharing systems can operate as a  system that has the advantages of both dockless and dock-based bikesharing systems: (1) it can meet the return demand of riders as in the case of dockless bikesharing systems and (2) it has the advantage of relatively easy bike management as in the case of dock-based bikesharing systems. To further minimize the cost, the capacity at each site can be set differently according to the demand. Therefore, the extension model can be used as an effective tool for designing the capacity in dock-based bikesharing systems.

\section{Conclusion and Discussions}\label{conclusion}
In this study, we defined a mixed discrete and continuous optimization problem and developed a network flow model to manage the daily operations of dockless bikesharing systems. We used the bikesharing trip data from Yishun, Singapore, as the input to perform a numerical experiment using the network flow model. Sensitivity analyses were then performed on the infrastructure scale, relocation cost, and usage price. We also considered both the optimal revenue and riders’ satisfaction level, which are both crucial indicators of the operation quality of a bikesharing system. As an extension to our model, we added a capacity constraint to each site to convert the dockless bikesharing relocation problem to a dock-based bikesharing relocation problem. Finally, we analyzed the relationship between dockless and dock-based bikesharing systems with this extended model.

For a system with a known number of bikes, the influence of trikes can be easily calculated using our network flow model to help operators manage the relocation capacity. However, this may be more complex when deciding the number of bikes in a brand new system. As illustrated in the previous sections, there are two ways for operators to increase the bike-riding demand. The first is to have more bikes in the system, and the other is to use more trikes for relocation. Since dockless bike sharing is an emerging service, to increase their market share, many operators tend to concentrate on market competition to place as many bikes as possible on the market, which is beyond necessary. The reason why some companies went bankrupt can be partly explained by this malignant competition; an oversupply of bikes leads to a much greater cost than the profit. Apart from causing businesses to collapse, the oversupply of bikes also requires too much public space, which makes city management harder. Although bikesharing has the advantage of reducing carbon emissions, the management of public resources is also considered an environmental issue. If the reduction of carbon emissions comes at the cost of wasting public space, then this service may not be a reliable option.

The extended model shows an intuitive result of the transformation from a dock-based bikesharing system to a dockless one by increasing the capacity of each site, which is achieved in reality by constructing more docks at each site. However, the number of docks needs to be carefully calculated and decided to control the construction cost.
 
From the analysis in this paper, we can make some policy suggestions for the management and operation of bikesharing systems. First, it is of practical implication for the operator to survey on the price influence before choosing its pricing strategy. There does not always exist a persistent positive correlation between the service price and revenue as the demand pattern is heavily dependent on the price. Thus, it is worthy of a comprehensive analysis of this impact for the operator to achieve higher revenue. What's more, higher revenue does not indicate a higher productivity of the system. Bikesharing is a business for the operator but also an important practice for guaranteeing social welfare. Therefore,  the customer satisfaction rate also needs to be paid special attention to. A best profit strategy does not indicate a higher service quality, so government oversight is indispensable to regulate the operator's operation. Second, a limitation of the number of bikes in the system needs to be imposed. The results obtained from the extended model have shown that reasonably setting the capacity limits at different sites does not seriously harm the revenue of the system. An oversupply of shared bikes in a dockless bikesharing system not only increases the operation costs but also wastes too much public space and makes city management harder. Therefore, governments should set limits on the total number of shared bikes allowed to be placed on the market. Doing this will not decrease the profits of bikesharing operators or reduce the service level.

As bikesharing systems are continuously developed in different cities, relocation problems are becoming more and more diverse and complex in different scenarios. Further research is needed to meet the demand for bikesharing relocation development. Specifically, the following aspects are proposed for future work. In our study, some assumptions were made for simplification to make the computation easier. In future research, more practical factors should be considered in the model, such as allowing waiting or moving to a nearby site when a rider arrives and finds no bike available, which has an impact on the satisfaction rate. This satisfaction rate also needs to be more carefully defined. In addition to using trikes as a relocation method, riding with a bonus is another effective strategy to consider for bike fleet relocation. Moreover, future research should focus on exploring more efficient algorithms that allow for an increase in the size of the problem. The problem size for a dockless bikesharing system can be much larger than a dock-based one because bikes can be returned anywhere in a dockless bikesharing system, requiring considering every place of the studied area. Finally, in the network flow model proposed in this work, the movement of each truck trike can be seen as a vehicle routing problem (VRP). Thus, future studies may concentrate on a comparative analysis between the network flow model and the VRP in a large-scale bikesharing system.

\section{Acknowledgements}
The research is sponsored by the Natural Science Foundation of China (Grant Nos. 71901164, 71777150, 71890970, 72021002), the Natural Science Foundation of Shanghai (Grant No. 19ZR1460700), and National Key Research and Development Program of China (Grant No. 2018YFB1600900). The corresponding author would like to acknowledge the support from Shanghai Pujiang Program (Grant No. 2019PJC107).





\bibliography{mylib.bib}







\end{document}